\documentclass[12pt,a4paper]{article}
\usepackage{}
\usepackage{amsfonts}
\usepackage{amssymb}
\usepackage{mathrsfs}
\usepackage{amsmath}
\usepackage{amsthm}
\usepackage{enumerate}
\usepackage{cite}
\usepackage[all]{xy}
\allowdisplaybreaks
\newtheorem{defn}{Definition}[section]
\newtheorem{thm}[defn]{Theorem}
\newtheorem{lem}[defn]{Lemma}
\newtheorem{prop}[defn]{Proposition}
\newtheorem{cor}[defn]{Corollary}
\newtheorem{eg}[defn]{Example}
\newtheorem{re}[defn]{Remark}

\newcommand\relphantom[1]{\mathrel{\phantom{#1}}}

\newcommand{\bdefn}{\begin{defn}}
\newcommand{\edefn}{\end{defn}}
\newcommand{\bthm}{\begin{thm}}
\newcommand{\ethm}{\end{thm}}
\newcommand{\blem}{\begin{lem}}
\newcommand{\elem}{\end{lem}}
\newcommand{\bprop}{\begin{prop}}
\newcommand{\eprop}{\end{prop}}
\newcommand{\bcor}{\begin{cor}}
\newcommand{\ecor}{\end{cor}}
\newcommand{\beg}{\begin{eg}}
\newcommand{\eeg}{\end{eg}}
\newcommand{\bre}{\begin{re}}
\newcommand{\ere}{\end{re}}
\newcommand{\bpf}{\begin{proof}}
\newcommand{\epf}{\end{proof}}

\newcommand{\K}{\mathbb{K}}

\newcommand{\Z}{\mathbb{Z}}
\newcommand{\g}{\mathfrak{g}}

\newcommand{\x}{\mathscr{X}}

\newcommand{\benu}{\begin{enumerate}}
\newcommand{\eenu}{\end{enumerate}}
\newcommand{\bc}{\begin{center}}
\newcommand{\ec}{\end{center}}
\newcommand{\bea}{\begin{eqnarray}}
\newcommand{\eea}{\end{eqnarray}}
\newcommand{\Bea}{\begin{eqnarray*}}
\newcommand{\Eea}{\end{eqnarray*}}
\newcommand{\beq}{\begin{equation}}
\newcommand{\eeq}{\end{equation}}
\newcommand{\Beq}{\begin{equation*}}
\newcommand{\Eeq}{\end{equation*}}
\newcommand{\bspl}{\begin{split}}
\newcommand{\espl}{\end{split}}

\numberwithin{equation}{section}
\newcommand{\supercite}[1]{\textsuperscript{\cite{#1}}}

\begin{document}
\title{{\bf  On the deformations and derivations of $n$-ary multiplicative Hom-Nambu-Lie superalgebras}}
\author{ Baoling Guan$^{1,2},$  Liangyun Chen$^{1}$
 \date{{\small {$^1$ School of Mathematics and Statistics, Northeast Normal
 University,\\
Changchun 130024, China}\\{\small {$^2$ College of Sciences, Qiqihar
University, Qiqihar 161006, China}}}}}

\maketitle
\date{}

\begin{abstract}

 In this paper, we introduce the relevant
 concepts of $n$-ary multiplicative Hom-Nambu-Lie superalgebras and construct three classes of $n$-ary multiplicative Hom-Nambu-Lie superalgebras. As a generalization of the notion of derivations for $n$-ary multiplicative Hom-Nambu-Lie algebras, we discuss the derivations of $n$-ary multiplicative Hom-Nambu-Lie superalgebras. In addition, the theory of one parameter formal deformation of $n$-ary multiplicative Hom-Nambu-Lie superalgebras is developed by choosing a suitable cohomology.
\bigskip

\end{abstract}
\renewcommand{\thefootnote}{\fnsymbol{footnote}}
\footnote[0]{ Corresponding author(L. Chen): chenly640@nenu.edu.cn.}
\footnote[0]{Supported by  NNSF of China (No.11171055),  NSF of  Jilin province (No.201115006), Scientific
Research Foundation for Returned Scholars
    Ministry of Education of China  and the Fundamental Research Funds for the Central Universities (No.12SSXT139). }

\section{Introduction}

In 1996, the concept of $n$-Lie
superalgebras was firstly introduced by Y. Daletskii and V. Kushnirevich in {\rm
  \cite{dykv}}. Moreover,  N. Cantarini and V. G. Kac gave a more general concept of $n$-Lie
superalgebras again  in 2010 in {\cite{CK}}. $n$-Lie superalgebras are more general structures
including $n$-Lie algebras ($n$-ary Nambu-Lie algebras), $n$-ary Nambu-Lie superalgebras and Lie superalgebras.

 The general Hom-algebra structures arose first in connection with quasi-deformation and discretizations of Lie algebras
of vector fields. These quasi-deformations lead to quasi-Lie algebras, a generalized Lie algebra structure in which the skewsymmetry
and Jacobi conditions are twisted. Hom-Lie algebras, Hom-associative algebras, Hom-Lie superalgebras, Hom-bialgebras, $n$-ary Hom-Nambu-Lie algebras and quasi-Hom-Lie algebras are discussed in \cite{aaf,ld,bs,masss,FSA,hls,af,mas,ydd,yddd,sy}. Generalizations of $n$-ary algebras of Lie type and associative type by twisting the identities using linear maps have been introduced in \cite{ah}.

The mathematical theory of deformations has proved to be a powerful tool in modeling physical reality. For example, (algebras associated with) classical quantum mechanics (and field theory) on a Poisson phase space can be deformed to (algebras associated with) quantum mechanics (and quantum field theory). The deformation of algebraic systems has been one of the problems that many mathematical researchers are interested in, Gerstenhaber studied the deformation theory of algebras in a series of papers \supercite{gm1,gm2,gm3,gm4,gm5}.  For example, it has been extended to covariant functors from a small category to algebras. In \cite{gmsds} and \cite{fmgmvaa},  it is respectively extended to algebra systems, bialgebras, Hopf algebras, Leibniz pairs and Poisson algebras, etc. In \cite{gm2}, Gerstenhaber developed a theory of deformation of
associative and Lie algebras. His theory links cohomologies of these algebras
and the Gerstenhaber bracket giving ¡®obstructions¡¯ to deformations. Nijenhuis
and Richardson noticed strong similarities between Gerstenhaber¡¯s theory and the
deformations of complex analytic structures on compact manifolds \supercite{narrw1}. They
axiomatized the theory of deformations via the introduction of graded Lie algebras
\supercite{narrw2}. One such example was given by the theory of deformations
of homomorphisms \supercite{narrw3}. Inspired by these works, we study the deformation theory of $n$-ary multiplicative Hom-Nambu-Lie superalgebras in this paper. In addition, the paper also discusses derivations of $n$-ary multiplicative Hom-Nambu-Lie superalgebras as a generalization of the notions of derivations for $n$-ary multiplicative Hom-Nambu-Lie algebras.

This paper is organized as follows. In section 1, we introduce the relevant
 concepts of $n$-ary multiplicative Hom-Nambu-Lie superalgebras and construct three classes of $n$-ary multiplicative Hom-Nambu-Lie superalgebras.
In section 2, the notion of derivation
introduced for $n$-ary multiplicative Hom-Nambu-Lie algebras in \cite{FSA} is extended to $n$-ary multiplicative Hom-Nambu-Lie superalgebras. In section 3, the theory of deformations of $n$-ary
multiplicative Hom-Nambu-Lie superalgebras is developed by choosing a suitable cohomology.
%\supercite{CK}

\bdefn \supercite{afn}
An $n$-ary Nambu-Lie superalgebra is a pair $(\g, [\cdot,\cdots,\cdot])$ consisting of a $\Z_2$-graded vector space $\g=\g_{\bar{0}}\oplus \g_{\bar{1}}$ and a multilinear mapping $[\cdot,\cdots,\cdot]: \underbrace{\g\times\cdots\times \g}_{n}\rightarrow \g,$ satisfying
\begin{align*}|[x_1,\cdots, x_n]|=&|x_1|+\cdots+|x_n|,\\
[x_1,\cdots, x_i, x_{i+1},\cdots, x_n]=&-(-1)^{|x_i| |x_{i+1}|}[x_1,\cdots, x_{i+1}, x_i,\cdots, x_n],\\
\begin{split}
[x_1,\cdots, x_{n-1}, [y_1,\cdots, y_n]]=&\sum_{i=1}^n(-1)^{(|x_1|+\cdots+|x_{n-1}|)(|y_1|+\cdots+|y_{i-1}|)}\\
\cdot[y_1,\cdots,y_{i-1},&[x_1,\cdots, x_{n-1},y_i], y_{i+1},\cdots,y_n],
\end{split}\end{align*}
where $|x|\in\Z_2$ denotes the degree of a homogeneous element $x\in\g$.
\edefn

\bdefn
An $n$-ary Hom-Nambu-Lie superalgebra is a triple $(\g, [\cdot,\cdots,\cdot],\alpha)$ consisting of a $\Z_2$-graded vector space $\g=\g_{\bar{0}}\oplus \g_{\bar{1}},$ a multilinear mapping $[\cdot,\cdots,\cdot]: \underbrace{\g\times\cdots\times \g}_{n}\rightarrow \g$ and a family $\alpha=(\alpha_{i})_{1\leq i\leq n-1}$ of even linear maps $\alpha_{i}: \g\rightarrow \g,$ satisfying
\begin{align}|[x_1,\cdots, x_n]|=&|x_1|+\cdots+|x_n|,\label{eq:1}\\
[x_1,\cdots, x_i, x_{i+1},\cdots, x_n]=&-(-1)^{|x_i| |x_{i+1}|}[x_1,\cdots, x_{i+1}, x_i,\cdots, x_n],\label{eq:2}\\
\begin{split}
[\alpha_{1}(x_1),\cdots, \alpha_{n-1}(x_{n-1}), [y_1,\cdots, y_n]]=&\sum_{i=1}^n(-1)^{(|x_1|+\cdots+|x_{n-1}|)(|y_1|+\cdots+|y_{i-1}|)}\label{eq:3}\\
\cdot[\alpha_{1}(y_1),\cdots,\alpha_{i-1}(y_{i-1}),&[x_1,\cdots, x_{n-1},y_i],\alpha_{i}(y_{i+1}),\cdots,\alpha_{n-1}(y_n)],
\end{split}\end{align}
where $|x|\in\Z_2$ denotes the degree of a homogeneous element $x\in\g.$

An $n$-ary Hom-Nambu-Lie superalgebra $(\g, [\cdot,\cdots,\cdot],\alpha)$ is multiplicative, if
 $\alpha=\\(\alpha_{i})_{1\leq i\leq n-1}$ with $\alpha_{1}=\cdots=\alpha_{n-1}=\alpha$ and satisfying
$$\alpha[x_{1},\cdots,x_{n}]=[\alpha(x_{1}),\cdots,\alpha(x_{n})], \forall x_{1},x_{2},\cdots,x_{n}\in \g.$$
\edefn
If the $n$-ary Hom-Nambu-Lie superalgebra $(\g, [\cdot,\cdots,\cdot],\alpha)$ is multiplicative, then the equation (\ref{eq:3}) can be read:
\begin{align}
\begin{split}
[\alpha(x_1),\cdots, \alpha(x_{n-1}), [y_1,\cdots, y_n]]=&\sum_{i=1}^n(-1)^{(|x_1|+\cdots+|x_{n-1}|)(|y_1|+\cdots+|y_{i-1}|)}\label{eq:31}\\
\cdot[\alpha(y_1),\cdots,\alpha(y_{i-1}),&[x_1,\cdots, x_{n-1},y_i],\alpha(y_{i+1}),\cdots,\alpha(y_n)].
\end{split}\tag{${1.3}^{'}$}\end{align}

It is clear that $n$-ary Hom-Nambu-Lie algebras and Hom-Lie
superalgebras are particular cases of $n$-ary Hom-Nambu-Lie superalgebras. In the sequel, when the notation ``$|x|$'' appears,
it means that $x$ is a homogeneous element of degree $|x|$.

\bdefn
Let $(\g, [\cdot,\cdots,\cdot],\alpha)$ and $(\g^{'}, [\cdot,\cdots,\cdot]^{'},\alpha^{'})$ be two $n$-ary Hom-Nambu-Lie superalgebras, where $\alpha=(\alpha_{i})_{1\leq i\leq n-1}$ and $\alpha^{'}=(\alpha^{'}_{i})_{1\leq i\leq n-1}.$ A linear map $f: \g\rightarrow \g$ is an $n$-ary Hom-Nambu-Lie superalgebra morphism if it satisfies
$$f[x_{1},\cdots,x_{n}]=[f(x_{1}),\cdots,f(x_{n})]^{'},$$
$$f\circ \alpha_{i}=\alpha_{i}^{'}\circ f, \forall i=1,\cdots,n-1.$$
\edefn

\bthm
Let $(\g, [\cdot,\cdots,\cdot], \alpha)$ be an $n$-ary multiplicative Hom-Nambu-Lie superalgebra and let $\beta: \g\rightarrow \g$ be a morphism of $g$ such that $\beta\circ\alpha=\alpha\circ\beta.$
Then $(\g, \beta\circ[\cdot,\cdots,\cdot],\beta\circ\alpha)$ is an $n$-ary multiplicative Hom-Nambu-Lie superalgebra.
\ethm
\bpf
Put $[\cdot,\cdots,\cdot]_{\beta}:=\beta\circ[\cdot,\cdots,\cdot].$ Then
\begin{align*}
&(\beta\circ\alpha)[x_{1},\cdots,x_{n}]_{\beta}=(\beta\circ\alpha)(\beta[x_{1},\cdots,x_{n}])\\
=&\beta\circ(\alpha\circ\beta)[x_{1},\cdots,x_{n}]\\
=&\beta\circ(\beta\circ\alpha)[x_{1},\cdots,x_{n}]\\
=&\beta[\beta\circ\alpha(x_{1}),\cdots, \beta\circ\alpha(x_{n})]\\
=&[\beta\circ\alpha(x_{1}),\cdots, \beta\circ\alpha(x_{n})]_{\beta},
\end{align*}
i.e., $\beta\circ\alpha$ is a morphism of $g.$
Moreover, we have
\begin{align*}
&[\beta\circ\alpha(x_{1}),\cdots,\beta\circ\alpha(x_{n-1}),[y_{1},\cdots,y_{n}]_{\beta}]_{\beta}\\
=&\beta[\beta\circ\alpha(x_{1}),\cdots,\beta\circ\alpha(x_{n-1}),{\beta}[y_{1},\cdots,y_{n}]]\\
=&\beta^{2}[\alpha(x_{1}),\cdots,\alpha(x_{n-1}),[y_{1},\cdots,y_{n}]]\\
=&\beta^{2}\big(\sum_{i=1}^{n}(-1)^{(|x_1|+\cdots+|x_{n-1}|)(|y_1|+\cdots+|y_{i-1}|)}[\alpha(y_{1}),\cdots,[x_{1},\cdots,x_{n-1},y_{i}],\cdots,\alpha(y_{n})]\big)\\
=&\sum_{i=1}^{n}(-1)^{(|x_1|+\cdots+|x_{n-1}|)(|y_1|+\cdots+|y_{i-1}|)}\beta[\beta\circ\alpha(y_{1}),\cdots,
\beta[x_{1},\cdots,x_{n-1},y_{i}],\cdots,\beta\circ\alpha(y_{n})]\\
=&\sum_{i=1}^{n}(-1)^{(|x_1|+\cdots+|x_{n-1}|)(|y_1|+\cdots+|y_{i-1}|)}[\beta\circ\alpha(y_{1}),\cdots,
[x_{1},\cdots,x_{n-1},y_{i}]_{\beta},\cdots,\beta\circ\alpha(y_{n})]_{\beta}.
\end{align*}
Therefore, $(\g, \beta\circ[\cdot,\cdots,\cdot],\beta\circ\alpha)$ is an $n$-ary multiplicative Hom-Nambu-Lie superalgebra.
\epf
In particular, we have the following example.
\beg
Let $(\g, [\cdot,\cdots,\cdot])$ be an $n$-ary Nambu-Lie superalgebra and let $\rho: \g\rightarrow \g$ be  an $n$-ary Nambu-Lie superalgebra endomorphism.
Then $(\g, \rho\circ[\cdot,\cdots,\cdot],\rho)$ is an $n$-ary multiplicative Hom-Nambu-Lie superalgebra.
\eeg
\bdefn Let $(\g,[\cdot,\cdots,\cdot]_{\g},\alpha)$ be an $n$-ary Hom-Nambu-Lie superalgebra.  A graded subspace $H\subseteq \g$ is a Hom-subalgebra of $(\g,[\cdot,\cdots,\cdot]_{\g},\alpha)$ if $\alpha(H)\subseteq H$
and $H$ is closed under the bracket operation $[\cdot,\cdots,\cdot]_{\g},$ i.e., $[u_{1},\cdots,u_{n}]_{\g}\in H, \forall u_{1},\cdots,u_{n}\in H.$\\
A graded subspace $H\subseteq \g$ is a Hom-ideal of $(\g,[\cdot,\cdots,\cdot]_{\g},\alpha)$ if $\alpha(H)\subseteq H$
and $[u_{1},u_{2},\cdots,\\u_{n}]_{\g}\in H, \forall u_{1}\in H, u_{2},\cdots,u_{n}\in \g.$
\edefn

\bdefn Let $(\g_{1},[\cdot,\cdots,\cdot]_{1},\alpha)$ and $(\g_{2},[\cdot,\cdots,\cdot]_{2},\beta)$ be two $n$-ary multiplicative Hom-Nambu-Lie superalgebras. Suppose that $\phi:\g_{1} \rightarrow \g_{2}$ is a linear map. $\mathfrak{G}_\phi=\{(x,\phi(x))|x\in \g_{1}\}\subseteq \g_{1}\oplus\g_{2}$ is called as the graph of a linear map $\phi:\g_{1} \rightarrow \g_{2}.$
\edefn

\bprop\label{proposition2.1}
Given two $n$-ary multiplicative Hom-Nambu-Lie superalgebras $(\g_{1},[\cdot,\\\cdots,\cdot]_{\g_{1}},\alpha)$ and $(\g_{2},[\cdot,\cdots,\cdot]_{\g_{2}},\beta)$, there is an $n$-ary multiplicative Hom-Nambu-Lie superalgebra
$(\g_{1}\oplus\g_{2},[\cdot,\cdots,\cdot]_{\g_{1}\oplus\g_{2}},\alpha+\beta)$, where the bilinear map $[\cdot,\cdots,\cdot]_{\g_{1}\oplus\g_{2}}:{\wedge}^2(\g_{1}\oplus\g_{2})\rightarrow \g_{1}\oplus\g_{2}$ is given by
$${[u_1+v_1,\cdots,u_n+v_n]}_{\g_{1}\oplus\g_{2}}={[u_1,\cdots,u_n]}_{\g_{1}}+{[v_1,\cdots,v_n]}_{\g_{2}},\forall  u_i\in \g_{1},  v_i \in \g_{2} (i=1,2,\cdots,n) $$
and the linear map $(\alpha+\beta):\g_{1}\oplus\g_{2}\rightarrow \g_{1}\oplus\g_{2}$ is given by
$$(\alpha+\beta)(u+v)=\alpha(u)+\beta(v),\forall u\in \g_{1}, v\in \g_{2}.$$
\eprop
\bpf
 For any  $u_i\in \g_{1},  v_i\in \g_{2},$ we have
\begin{align*}
&{[u_1+v_1,\cdots,u_i+v_i,u_{i+1}+v_{i+1},\cdots,u_n+v_n]}_{\g_{1}\oplus\g_{2}}\\
=&{[u_1,\cdots,u_i,u_{i+1},\cdots,u_n]}_{\g_{1}}+{[v_1,\cdots,v_i,v_{i+1},\cdots,v_n]}_{\g_{2}}\\
=&-(-1)^{|u_i||u_{i+1}|}{[u_1,\cdots,u_{i+1},u_i,\cdots,u_n]}_{\g_{1}}
-(-1)^{|v_i||v_{i+1}|}{[v_1,\cdots,v_{i+1},v_i,\cdots,v_n]}_{\g_{2}}\\
=&-(-1)^{|u_i||u_{i+1}|}({[u_1,\cdots,u_i,u_{i+1},\cdots,u_n]}_{\g_{1}}
+{[v_1,\cdots,v_i,v_{i+1},\cdots,v_n]}_{\g_{2}})\\
=&-(-1)^{|u_i||u_{i+1}|}{[u_1+v_1,\cdots,u_{i+1}+v_{i+1},u_i+v_i,\cdots,u_n+v_n]}_{\g_{1}\oplus\g_{2}}.
\end{align*}
The bracket is obviously supersymmetric. By a direct computation we
have
\begin{align*}
&[(\alpha+\beta)(u_1+v_1),\cdots, (\alpha+\beta)(u_{n-1}+v_{n-1}), [x_{1}+y_1,\cdots, x_{n}+y_n]_{\g_{1}\oplus\g_{2}}]_{\g_{1}\oplus\g_{2}}\\
=&[\alpha(u_1)+\beta(v_1),\cdots, \alpha(u_{n-1})+\beta(v_{n-1}), [x_{1}+y_1,\cdots, x_{n}+y_n]_{\g_{1}\oplus\g_{2}}]_{\g_{1}\oplus\g_{2}}\\
=&[\alpha(u_1),\cdots, \alpha(u_{n-1}), [x_{1},\cdots, x_{n}]_{\g_{1}}]_{\g_{1}}+[\beta(v_1),\cdots, \beta(v_{n-1}), [y_1,\cdots, y_n]_{\g_{2}}]_{\g_{2}}\\
=&\sum_{i=1}^n(-1)^{(|u_1|+\cdots+|u_{n-1}|)(|x_1|+\cdots+|x_{i-1}|)}([\alpha(x_{1}),\cdots,\alpha(x_{i-1}),[u_1,\cdots, u_{n-1},x_i]_{\g_{1}},\\&\alpha(x_{i+1}),\cdots,\alpha(x_n)]_{\g_{1}}
+[\beta(y_1),\cdots,\beta(y_{i-1}),[v_1,\cdots, v_{n-1},y_i]_{\g_{2}},\beta(y_{i+1}),\cdots,\beta(y_n)]_{\g_{2}})\\
=&\sum_{i=1}^n(-1)^{(|u_1|+\cdots+|u_{n-1}|)(|x_1|+\cdots+|x_{i-1}|)}[\alpha(x_{1})+\beta(y_1),\cdots,\alpha(x_{i-1})+\beta(y_{i-1}),\\
&[u_1,\cdots, u_{n-1},x_i]_{\g_{1}}+[v_1,\cdots, v_{n-1},y_i]_{\g_{2}},\alpha(x_{i+1})+\beta(y_{i+1}),\cdots,\alpha(x_n)+\beta(y_n)]_{\g_{1}\oplus\g_{2}}\\
=&\sum_{i=1}^n(-1)^{(|u_1|+\cdots+|u_{n-1}|)(|x_1|+\cdots+|x_{i-1}|)}[(\alpha+\beta)(x_{1}+y_1),\cdots,(\alpha+\beta)(x_{i-1}+y_{i-1}),\\
&[u_1+v_1,\cdots, u_{n-1}+v_{n-1},x_i+y_i]_{\g_{1}\oplus\g_{2}},(\alpha+\beta)(x_{i+1}+y_{i+1}),\cdots,(\alpha+\beta)(x_n+y_n)]_{\g_{1}\oplus\g_{2}}\\
=&\sum_{i=1}^n(-1)^{(|u_1|+\cdots+|u_{n-1}|)(|x_1|+\cdots+|x_{i-1}|)}[(\alpha+\beta)(x_{1}+y_1),\cdots,(\alpha+\beta)(x_{i-1}+y_{i-1}),\\
&[u_1+v_1,\cdots, u_{n-1}+v_{n-1},x_i+y_i]_{\g_{1}\oplus\g_{2}},(\alpha+\beta)(x_{i+1}+y_{i+1}),\cdots,(\alpha+\beta)(x_n+y_n)]_{\g_{1}\oplus\g_{2}}.
\end{align*}
\epf

\bprop\label{proposition2.1}
A linear map $\phi:(\g_{1},[\cdot,\cdots,\cdot]_{\g_{1}},\alpha)\rightarrow(\g_{2},[\cdot,\cdots,\cdot]_{\g_{2}},\beta)$ is a morphism of $n$-ary multiplicative Hom-Nambu-Lie superalgebras if and only if the graph $\mathfrak{G}_\phi\subseteq \g_{1}\oplus\g_{2}$ is a Hom-subalgebra of  $(\g_{1}\oplus\g_{2},[\cdot,\cdots,\cdot]_{\g_{1}\oplus\g_{2}},\alpha+\beta)$.
\eprop
\bpf
Let
$\phi:(\g_{1},[\cdot,\cdots,\cdot]_{\g_{1}},\alpha)\rightarrow(\g_{2},[\cdot,\cdots,\cdot]_{\g_{2}},\beta)$
be a morphism of $n$-ary multiplicative Hom-Nambu-Lie superalgebras. Then
\begin{align*}
&[u_{1}+\phi(u_{1}),\cdots,u_{n}+\phi(u_{n})]_{\g_{1}\oplus\g_{2}}\\
=&[u_{1},\cdots,u_{n}]_{\g_{1}}+[\phi(u_{1}),\cdots,\phi(u_{n})]_{\g_{2}}\\
=&[u_{1},\cdots,u_{n}]_{\g_{1}}+\phi[u_{1},\cdots,u_{n}]_{\g_{2}}.
\end{align*}

Then the graph $\mathfrak{G}_\phi$ is closed under the bracket operation $[\cdot,\cdots,\cdot]_{\g_{1}\oplus\g_{2}}$. Furthermore, we obtain
$$(\alpha+\beta)(u+\phi(u))=\alpha(u)+\beta\circ\phi(u)=\alpha(u)+\phi\circ\alpha(u),$$
which implies that $(\alpha+\beta)(\mathfrak{G}_\phi)\subseteq \mathfrak{G}_\phi$. Thus, $\mathfrak{G}_\phi$ is a Hom-subalgebra of $(\g_{1}\oplus\g_{2},[\cdot,\cdots,\cdot]_{\g_{1}\oplus\g_{2}},\\ \alpha+\beta)$.

Conversely, if the graph $\mathfrak{G}_\phi\subseteq \g_{1}\oplus\g_{2}$ is a Hom-subalgebra of
$(\g_{1}\oplus\g_{2},[\cdot,\cdots,\cdot]_{\g_{1}\oplus\g_{2}},\alpha+\beta)$, then we have
$$[u_{1}+\phi(u_{1}),\cdots,u_{n}+\phi(u_{n})]_{\g_{1}\oplus\g_{2}}=
[u_{1},\cdots,u_{n}]_{\g_{1}}+\phi[u_{1},\cdots,u_{n}]_{\g_{2}},$$
which implies that $$\phi[u_{1},\cdots,u_{n}]_{\g_{2}}=\phi[u_{1},\cdots,u_{n}]_{\g_{1}}.$$
Furthermore, $(\alpha+\beta)(\mathfrak{G}_{\phi})\subset \mathfrak{G}_{\phi}$ yields that
$$(\alpha+\beta)(u+\phi(u))=\alpha(u)+\beta\circ\phi(u)\in\mathfrak{G}_\phi,$$
which is equivalent to the condition $\beta\circ\phi(u)=\phi\circ\alpha(u)$, i.e. $\beta\circ\phi=\phi\circ\alpha.$ Therefore, $\phi$ is a morphism of $n$-ary multiplicative Hom-Nambu-Lie superalgebras.
\epf

\section{Derivations of $n$-ary multiplicative Hom-Nambu-Lie superalgebras}

Let $(\g, [\cdot,\cdots,\cdot],\alpha)$ be an $n$-ary multiplicative Hom-Nambu-Lie superalgebra. We denote by $\alpha^{k}$ the $k$-times compositions of $\alpha.$
In particular, we set $\alpha^{0}=\mathrm{id}.$
\bdefn
For $k\geq 0,$ we call $D\in \mathrm{End}(\g)$ an $\alpha^{k}$-derivation of the $n$-ary multiplicative Hom-Nambu-Lie superalgebra $(\g, [\cdot,\cdots,\cdot],\alpha)$
if $$D\circ\alpha=\alpha\circ D$$
and for $x_{i}\in \g(i=1,\cdots,n),$
$$D[x_{1},\cdots,x_{n}]=\sum_{i=1}^n(-1)^{|D|(|x_{1}|+\cdots+|x_{i-1}|)}[\alpha^{k}(x_{1}),\cdots,\alpha^{k}(x_{i-1}),D(x_{i}),\alpha^{k}(x_{i+1}),\cdots,\alpha^{k}(x_{n})].$$
\edefn
We denote by $\mathrm{Der}_{\alpha^{k}}(\g)$ the set of $\alpha^{k}$-derivations of the $n$-ary multiplicative Hom-Nambu-Lie superalgebra $(\g, [\cdot,\cdots,\cdot],\alpha).$ Notice that we obtain classical derivations for $k=0.$

For $\x\in \g^{\wedge^{n-1}}$ satisfying $\alpha(\x)=\x$ and $k\geq 0,$ we define the map $\mathrm{ad}_{k}(\x)\in \mathrm{End}(\g)$ by
$$\mathrm{ad}_{k}(\x)(y)=[x_{1},\cdots,x_{n-1},\alpha^{k}(y)], \forall y\in \g.$$
Then

\blem
The map $\mathrm{ad}_{k}(\x)$ is an $\alpha^{k+1}$-derivation and is called an inner $\alpha^{k+1}$-derivation.
\elem
We denote by $\mathrm{Inn}_{\alpha^{k}}(\g)$ the $\K$-vector space generated by all inner $\alpha^{k+1}$-derivations.
For any $D\in \mathrm{Der}_{\alpha^{k}}(\g)$ and $D^{'}\in \mathrm{Der}_{\alpha^{k^{'}}}(\g),$ we define their commutator $[D,D^{'}]=D\circ D^{'}-(-1)^{|D||D^{'}|}D^{'}\circ D.$ Set $\mathrm{Der}(\g)=\oplus_{k\geq 0}\mathrm{Der}_{\alpha^{k}}(\g)$ and $\mathrm{Inn}(\g)=\bigoplus_{k\geq 0}\mathrm{Inn}_{\alpha^{k}}(\g).$

\blem \label{5.3}
For any $D\in \mathrm{Der}_{\alpha^{k}}(\g)$ and $D^{'}\in \mathrm{Der}_{\alpha^{k^{'}}}(\g),$ we have $[D,D^{'}]\in \mathrm{Der}_{\alpha^{k+k^{'}}}(\g).$
\elem
\bpf
Let $x_{i}\in \g, 1\leq i\leq n.$ $D\in \mathrm{Der}_{\alpha^{k}}(\g)$ and $D^{'}\in \mathrm{Der}_{\alpha^{k^{'}}}(\g),$ then
\begin{align*}
&D\circ D^{'}([x_{1},\cdots,x_{n}])\\
=&D(\sum_{i=1}^n(-1)^{|D^{'}|(|x_{1}|+\cdots+|x_{i-1}|)}
[\alpha^{k^{'}}(x_{1}),\cdots,\alpha^{k^{'}}(x_{i-1}),D^{'}(x_{i}),\alpha^{k^{'}}(x_{i+1}),\cdots,\alpha^{k^{'}}(x_{n})])\\
=&\sum_{i=1}^n(-1)^{|D^{'}|(|x_{1}|+\cdots+|x_{i-1}|)}
D[\alpha^{k^{'}}(x_{1}),\cdots,\alpha^{k^{'}}(x_{i-1}),D^{'}(x_{i}),\alpha^{k^{'}}(x_{i+1}),\cdots,\alpha^{k^{'}}(x_{n})]\\
 =&\sum_{i=1}^n(-1)^{|D^{'}|(|x_{1}|+\cdots+|x_{i-1}|)}(-1)^{|D|(|x_{1}|+\cdots+|x_{i-1}|)}\\
 &\cdot [\alpha^{k+k^{'}}(x_{1}),\cdots,\alpha^{k+k^{'}}(x_{i-1}),D\circ D^{'}(x_{i}),\alpha^{k+k^{'}}(x_{i+1}),\cdots,\alpha^{k+k^{'}}(x_{n})]\\
 +&\sum_{i<j}(-1)^{|D^{'}|(|x_{1}|+\cdots+|x_{i-1}|)}(-1)^{|D|(|x_{1}|+\cdots+|x_{j-1}|+|D^{'}|)}\\
 &\cdot [\alpha^{k+k^{'}}(x_{1}),\cdots,\alpha^{k}(D^{'}(x_{i})),\cdots,\alpha^{k^{'}}(D(x_{j})),\cdots,\alpha^{k+k^{'}}(x_{n})]\\
 +&\sum_{i>j}(-1)^{|D^{'}|(|x_{1}|+\cdots+|x_{i-1}|)}(-1)^{|D|(|x_{1}|+\cdots+|x_{j-1}|)}\\
 &\cdot [\alpha^{k+k^{'}}(x_{1}),\cdots,\alpha^{k^{'}}(D(x_{j})),\cdots,\alpha^{k}(D^{'}(x_{i})),\cdots,\alpha^{k+k^{'}}(x_{n})].
\end{align*}
and
\begin{align*}
&-(-1)^{|D||D^{'}|}D^{'}\circ D([x_{1},\cdots,x_{n}])\\
=&-(-1)^{|D||D^{'}|}D^{'}(\sum_{i=1}^n(-1)^{|D|(|x_{1}|+\cdots+|x_{i-1}|)}
[\alpha^{k}(x_{1}),\cdots,D(x_{i}),\cdots,\alpha^{k}(x_{n})])\\
=&-(-1)^{|D||D^{'}|}\sum_{i=1}^n(-1)^{|D|(|x_{1}|+\cdots+|x_{i-1}|)}
D^{'}[\alpha^{k}(x_{1}),\cdots,D(x_{i}),\cdots,\alpha^{k}(x_{n})]\\
 =&-(-1)^{|D||D^{'}|}\sum_{i=1}^n(-1)^{|D|(|x_{1}|+\cdots+|x_{i-1}|)}(-1)^{|D^{'}|(|x_{1}|+\cdots+|x_{i-1}|)}\\
 &\cdot [\alpha^{k^{'}+k}(x_{1}),\cdots,\alpha^{k^{'}+k}(x_{i-1}),D^{'}\circ D(x_{i}),\alpha^{k^{'}+k}(x_{i+1}),\cdots,\alpha^{k^{'}+k}(x_{n})]\\
-&(-1)^{|D||D^{'}|}\sum_{i<j}(-1)^{|D|(|x_{1}|+\cdots+|x_{i-1}|)}(-1)^{|D^{'}|(|x_{1}|+\cdots+|x_{j-1}|+|D|)}\\
 &\cdot [\alpha^{k^{'}+k}(x_{1}),\cdots,\alpha^{k^{'}}(D(x_{i})),\cdots,\alpha^{k}(D^{'}(x_{j})),\cdots,\alpha^{k^{'}+k}(x_{n})]\\
 -&(-1)^{|D||D^{'}|}\sum_{i>j}(-1)^{|D|(|x_{1}|+\cdots+|x_{i-1}|)}(-1)^{|D^{'}|(|x_{1}|+\cdots+|x_{j-1}|)}\\
 &\cdot [\alpha^{k^{'}+k}(x_{1}),\cdots,\alpha^{k}(D^{'}(x_{j})),\cdots,\alpha^{k^{'}}(D(x_{i})),\cdots,\alpha^{k^{'}+k}(x_{n})].
\end{align*}
Then we obtain
\begin{align*}
&[D,D^{'}]([x_{1},\cdots,x_{n}])=(D\circ D^{'}-(-1)^{|D||D^{'}|}D^{'}\circ D)([x_{1},\cdots,x_{n}])\\
=&\sum_{i=1}^n(-1)^{|D^{'}|(|x_{1}|+\cdots+|x_{i-1}|)}(-1)^{|D|(|x_{1}|+\cdots+|x_{i-1}|)}[\alpha^{k+k^{'}}(x_{1}),\cdots,\alpha^{k+k^{'}}(x_{i-1}),\\
 &(D\circ D^{'}-(-1)^{|D||D^{'}|}D^{'}\circ D)(x_{i}),\alpha^{k+k^{'}}(x_{i+1}),\cdots,\alpha^{k+k^{'}}(x_{n})]\\
=&\sum_{i=1}^n(-1)^{(|D^{'}|+|D|)(|x_{1}|+\cdots+|x_{i-1}|)}\\
&\cdot[\alpha^{k+k^{'}}(x_{1}),\cdots,\alpha^{k+k^{'}}(x_{i-1}),[D,D^{'}](x_{i}),\alpha^{k+k^{'}}(x_{i+1}),\cdots,\alpha^{k+k^{'}}(x_{n})],
\end{align*}
which yields that $[D,D^{'}]\in \mathrm{Der}_{\alpha^{k+k^{'}}}(\g).$
\epf
\bprop
The pair $(\mathrm{Der}(\g), [\cdot,\cdot]),$ where the bracket is the usual commutator, defines a Lie superalgebra and  $\mathrm{Inn}(\g)$
constitutes an ideal of it.
\eprop

\bpf $(\mathrm{Der}(\g), [\cdot,\cdot])$ is a  Lie superalgebra by using Lemma \ref{5.3}. We show that $\mathrm{Inn}(\g)$
is an ideal. Let $\mathrm{ad}_{k-1}(\x)(y)=[x_{1},\cdots,x_{n-1},\alpha^{k-1}(y)]$ be an inner $\alpha^{k}$-derivation on $\g$ and $D\in \mathrm{Der}_{\alpha^{k^{'}}}(\g)$ for $k\geq 1$ and $k^{'}\geq 0.$ Then
 $[D,\mathrm{ad}_{k-1}(\x)]\in \mathrm{Der}_{\alpha^{k+k^{'}}}(\g)$
and for any $y\in \g$
\begin{align*}
&[D,\mathrm{ad}_{k-1}(\x)](y)\\
=&D[x_{1},\cdots,x_{n-1},\alpha^{k-1}(y)]
-(-1)^{|D|(|x_{1}|+\cdots+|x_{n-1}|)}[x_{1},\cdots,x_{n-1},\alpha^{k-1}(D(y))]\\
=&D[\alpha^{k}(x_{1}),\cdots,\alpha^{k}(x_{n-1}),\alpha^{k-1}(y)]
-(-1)^{|D|(|x_{1}|+\cdots+|x_{n-1}|)}\\
&\cdot[\alpha^{k+k^{'}}(x_{1}),\cdots,\alpha^{k+k^{'}}(x_{n-1}),\alpha^{k-1}(D(y))]\\
 =&\sum_{i\leq n-1}(-1)^{|D|(|x_{1}|+\cdots+|x_{i-1}|)} [\alpha^{k+k^{'}}(x_{1}),\cdots,D(\alpha^{k}(x_{i})),\cdots,\alpha^{k+k^{'}}(x_{n-1}),\alpha^{k+k^{'}-1}(y)]\\
 =&\sum_{i\leq n-1}(-1)^{|D|(|x_{1}|+\cdots+|x_{i-1}|)} [x_{1},\cdots,D(x_{i}),\cdots,x_{n-1},\alpha^{k+k^{'}-1}(y)]\\
 =&\sum_{i\leq n-1}(-1)^{|D|(|x_{1}|+\cdots+|x_{i-1}|)} \mathrm{ad}_{k+k^{'}-1}(x_{1}\wedge\cdots\wedge D(x_{i})\wedge\cdots\wedge x_{n-1})(y).
\end{align*}
Therefore, $[D,\mathrm{ad}_{k-1}(\x)]\in \mathrm{Inn}_{\alpha^{k+k^{'}-1}}(\g).$
\epf
\section{Deformations of $n$-ary multiplicative Hom-Nambu-Lie superalgebras}
\bdefn \supercite{gcm} For $m\geq 1,$ we call $m$-coboundary operator of the $n$-ary multiplicative Hom-Nambu-Lie superalgebra $(\g, [\cdot,\cdots,\cdot],\alpha)$
the even linear map $\delta^{m}: C^m(\g, V)\rightarrow C^{m+1}(\g, V)$ by
\begin{align*}
&(\delta^{m} f)(\x_1,\cdots,\x_m, \x_{m+1}, z)\\
=&\sum_{i<j}(-1)^i(-1)^{|\x_i|(|\x_{i+1}|+\cdots+|\x_{j-1}|)}f(\alpha(\x_1),\cdots,\widehat{\alpha(\x_i)},\cdots,[\x_{i},\x_{j}]_{\alpha},\cdots,\alpha(\x_{m+1}),\alpha(z))\\
&+\sum_{i=1}^{m+1}(-1)^i(-1)^{|\x_i|(|\x_{i+1}|+\cdots+|\x_{m+1}|)}f(\alpha(\x_1),\cdots,\widehat{\alpha(\x_i)},\cdots,\alpha(\x_{m+1}),\x_i\cdot z)\\
& +\sum_{i=1}^{m+1}(-1)^{i+1}(-1)^{|\x_i|(|f|+|\x_{1}|+\cdots+|\x_{i-1}|)}\alpha^{m}(\x_i)\cdot f(\x_1,\cdots,\widehat{\x_i},\cdots,\x_{m+1}, z)\\
&  +(-1)^m(f(\x_1,\cdots,\x_m, ~~)\cdot \x_{m+1})\bullet_{\alpha} \alpha^{m}(z),
\end{align*}
where $\x_i=\x_i^1\wedge\cdots\wedge\x_i^{n-1}\in\g^{\wedge^{n-1}}, i=1,\cdots,m+1, z\in\g$ and the last term is defined by
\Beq\begin{split}
(f(\x_1,\cdots,\x_m,~~ )\cdot \x_{m+1})\bullet_{\alpha} \alpha^{m}(z)=&\sum_{i=1}^{n-1}(-1)^{(|f|+|\x_{1}|+\cdots+|\x_{m}|)(|\x_{m+1}^1|+\cdots+|\x_{m+1}^{i-1}|)}\\
\cdot[\alpha^{m}(\x_{m+1}^1),\cdots,&f(\x_1,\cdots,\x_m,\x_{m+1}^i),\cdots,\alpha^{m}(\x_{m+1}^{n-1}),\alpha^{m}(z)].
\end{split}\Eeq
\edefn
\bthm \supercite{gcm}
Let $f\in C^{m}(\g, V)$ be an $m$-cochain. Then $\delta^{m+1}\circ \delta^{m}(f)=0.$
\ethm

In \cite{gcm}, it also points out that, the map $f\in C^{m}(\g, V)$ is called an $m$-supercocycle if $\delta^{m} f=0$. We denote by $Z^{m}(\g,V)$ the graded subspace spanned by $m$-supercocycles. Since $\delta^{m+1}\circ \delta^{m}(f)=0$ for all $f \in C^{m}(\g, V)$, $\delta^{m-1} C^{m-1}(\g, V)$ is a graded subspace of $Z^{m}(\g,V)$. Therefore we can define a graded cohomology space $H^{m}(\g,V)$ of $\g$ as the graded factor space $Z^{m}(\g,V)/\delta^{m-1} C^{m-1}(\g, V).$

We next will discuss the deformation of $n$-ary multiplicative Hom-Nambu-Lie superalgebras.
Let $\mathbb{K}[[t]]$ denote the power series ring in one variable $t$ with coefficients in $\mathbb{K}$ and $\g[[t]]$ be the set of
formal series whose coefficients are elements of the vector space $\g.$
\bdefn
Let $(\g, [\cdot,\cdots,\cdot],\alpha)$ be an $n$-ary multiplicative Hom-Nambu-Lie superalgebra over $\K$. A deformation of $(\g, [\cdot,\cdots,\cdot],\alpha)$ is given by $\K[[t]]$-$n$-linear map
$$f_{t}=\sum_{p\geq 0}f_{p}t^{p}: \g[[t]]\times\cdots\times\g[[t]]\rightarrow\g[[t]]$$
such that $(\g[[t]],f_{t},\alpha)$ is also an $n$-ary multiplicative Hom-Nambu-Lie superalgebra. We call $f_{1}$ the infinitesimal deformation of $(\g, [\cdot,\cdots,\cdot],\alpha).$
\edefn
Since $(\g[[t]],f_{t},\alpha)$ is an $n$-ary multiplicative Hom-Nambu-Lie superalgebra,
$f_{t}$ satisfies
\begin{align}\alpha  \circ f_{t}(x_{1},\cdots,x_{n})=&f_{t}(\alpha(x_{1}),\cdots,\alpha(x_{n})),\tag{a}\\
|f_{t}(x_{1},\cdots,x_{n})|=&|x_{1}|+\cdots+|x_{n}|,\tag{b}\\
\begin{split}f_{t}(\alpha(x_{1}),\cdots,\alpha(x_{n-1}),f_{t}(y_{1},\cdots,y_{n}))
=&\sum_{i=1}^n(-1)^{(|x_{1}|+\cdots+|x_{n-1}|)(|y_{1}|+\cdots+|y_{i-1}|)}\\
\quad\quad\quad\cdot f_{t}(\alpha(y_{1}),\cdots,\alpha(y_{i-1}),f_{t}(x_{1},&\cdots,x_{n-1},y_{i}),\alpha(y_{i+1}),\cdots,\alpha(y_{n})).
\end{split}\tag{c}
\end{align}
$\mathrm{(a)}-\mathrm{(c)}$ is respectively equivalent to
\begin{align}\alpha  \circ f_{p}(x_{1},\cdots,x_{n})=&f_{p}(\alpha(x_{1}),\cdots,\alpha(x_{n})),\tag{$a^{'}$}\\
|f_{p}(x_{1},\cdots,x_{n})|=&|x_{1}|+\cdots+|x_{n}|,\tag{$b^{'}$}\\
\begin{split}\sum_{p+q=l}f_{p}(\alpha(x_{1}),\cdots,\alpha(x_{n-1}),f_{q}(y_{1},\cdots,y_{n}))
=&\sum_{i=1}^n(-1)^{(|x_{1}|+\cdots+|x_{n-1}|)(|y_{1}|+\cdots+|y_{i-1}|)}\\
\quad\quad\quad\cdot(\sum_{p+q=l}f_{p}(\alpha(y_{1}),\cdots,\alpha(y_{i-1}),f_{q}(x_{1},&\cdots,x_{n-1},y_{i}),\alpha(y_{i+1}),\cdots,\alpha(y_{n})).
\end{split} \tag{$c^{'}$}
\end{align}
We call these the deformation equations for an $n$-ary multiplicative Hom-Nambu-Lie superalgebra.

($a^{'}$) and ($b^{'}$) shows that $f_{p}\in C^{1}(\g,\g)_{\bar{0}}.$ In ($c^{'}$), set $l=1,$ then
\begin{align*}
&[\alpha(x_{1}),\cdots,\alpha(x_{n-1}),f_{1}(y_{1},\cdots,y_{n})]+f_{1}(\alpha(x_{1}),\cdots,\alpha(x_{n-1}),[y_{1},\cdots,y_{n}])\\
&-\sum_{i=1}^n(-1)^{(|x_{1}|+\cdots+|x_{n-1}|)(|y_{1}|+\cdots+|y_{i-1}|)}\Big([\alpha(y_{1}),\cdots,\alpha(y_{i-1}),
f_{1}(x_{1},\cdots,x_{n-1},y_{i}),\\
&\quad \alpha(y_{i+1}),\cdots,\alpha(y_{n})]
+f_{1}(\alpha(y_{1}),\cdots,\alpha(y_{i-1}),[x_{1},\cdots,x_{n-1},y_{i}],\alpha(y_{i+1}),\cdots,\alpha(y_{n}))\Big)\\
&=0,
\end{align*}
i.e., $\delta^{1}f_{1}(x_{1},\cdots,x_{n-1},[y_{1},\cdots,y_{n-1},y_{n}])=0.$ Hence the infinitesimal deformation $f_{1}\in Z^{1}(\g,\g)_{\bar{0}}.$

\bdefn
Two deformations $f_{t}$ and $f_{t^{'}}$ of the $n$-ary multiplicative Hom-Nambu-Lie superalgebra $(\g, [\cdot,\cdots,\cdot],\alpha)$ are said to be equivalent, if there exists an isomorphism of $n$-ary multiplicative Hom-Nambu-Lie superalgebras $\Phi_{t}: (\g, f_{t},\alpha)\rightarrow (\g, f^{'}_{t},\alpha),$ where $\Phi_{t}=\sum_{i\geq0}\varphi_{i}t^{i}, \varphi_{i}: \g\rightarrow \g$ is a linear map such that $$\varphi_{0}=\mathrm{id}_{\g};\ \varphi_{i}\circ \alpha=\alpha\circ\varphi_{i};$$$$ \Phi_{t}\circ f_{t}(x_{1},\cdots,x_{n})=f^{'}_{t}(\Phi_{t}(x_{1}),\cdots,\Phi_{t}(x_{n})),$$ and is denoted by $f_{t}\sim f_{t^{'}}.$ When $f_{1}=f_{2}=\cdots=0,$ $f_{1}=f_{0}$ is called the null deformation; if $f_{t}\sim f_{0},$ then $f_{t}$ is called  the trivial deformation.
\edefn

\bthm
Let $f_{t}$ and $f_{t^{'}}$ be two equivalent deformations of the $n$-ary multiplicative Hom-Nambu-Lie superalgebra $(\g, [\cdot,\cdots,\cdot],\alpha).$
Then the infinitesimal deformations $f_{1}$ and $f^{'}_{1}$ belong to the same cohomology class in the cohomology group $H^{2}(\g,\g).$
\ethm
\bpf Put $B^{2}(\g,\g):=\delta^{1}C^{1}(\g,\g).$ It is enough to prove that $f_{1}-f^{'}_{1}\in B^{2}(\g,\g).$
Let $\Phi_{t}: (\g, f_{t},\alpha)\longrightarrow (\g, f^{'}_{t},\alpha)$ be an isomorphism of $n$-ary multiplicative Hom-Nambu-Lie superalgebras. Then $\varphi_{1}\in C^{1}(\g,\g)_{\bar{0}}$ and $$\sum_{i\geq0}\varphi_{i}(\sum_{j\geq0}f_{j}(x_{1},\cdots,x_{n}))t^{i+j}=\sum_{i\geq0}f^{'}_{i}(\sum_{j_{1}\geq0}\varphi_{j_{1}}(x_{1}),
\cdots,\sum_{j_{n}\geq0}\varphi_{j_{n}}(x_{n}))t^{i+j_{1}+\cdots+j_{n}},$$
comparing with the coefficients of $t^{1}$ for two sides of the above equation, we obtain
\begin{align*}
f_{1}(x_{1},\cdots,x_{n})+\varphi_{1}[x_{1},x_{2},\cdots,x_{n}]=&[\varphi_{1}(x_{1}),x_{2},\cdots,x_{n}]+
[x_{1},\varphi_{1}(x_{2}),x_{3},\cdots,x_{n}]\\
+&\cdots+[x_{1},\cdots,x_{n-1},\varphi_{1}(x_{n})]+f^{'}_{1}(x_{1},\cdots,x_{n}).
\end{align*}
Furthermore, one gets
\begin{align*}
&f_{1}(x_{1},\cdots,x_{n})-f^{'}_{1}(x_{1},\cdots,x_{n})=-\varphi_{1}[x_{1},x_{2},\cdots,x_{n}]+[\varphi_{1}(x_{1}),x_{2},\cdots,x_{n}]\\
&\qquad \qquad\qquad \qquad\quad+[x_{1},\varphi_{1}(x_{2}),x_{3},\cdots,x_{n}]+\cdots+[x_{1},\cdots,x_{n-1},\varphi_{1}(x_{n})]\\
&=-\varphi_{1}[x_{1},\cdots,x_{n}]+\sum_{i=1}^n(-1)^{n-i}(-1)^{|x_{i}|(|x_{i+1}|+\cdots+|x_{n}|)}[x_{1},\cdots,\widehat{x_{i}},\cdots,x_{n},\varphi_{1}(x_{i})]\\
&=\delta^{1}\varphi_{1}(x_{1},\cdots,x_{n}).
\end{align*}
Therefore, $f_{1}-f^{'}_{1}=\delta^{1}\varphi_{1}\in \delta^{1}C^{1}(\g,\g)_{\bar{0}},$ i.e., $f_{1}-f^{'}_{1}\in B^{2}(\g,\g).$
\epf

An $n$-ary multiplicative Hom-Nambu-Lie superalgebra $(\g, [\cdot,\cdots,\cdot],\alpha)$ is {\it analytically rigid} if every deformation $f_t$ is equivalent to the null deformation $f_0$. We have a fundamental theorem.
\bthm
If $(\g, [\cdot,\cdots,\cdot],\alpha)$ is an $n$-ary multiplicative Hom-Nambu-Lie superalgebra with $H^2(\g, \g)=0,$ then  $(\g, [\cdot,\cdots,\cdot],\alpha)$ is analytically rigid.
\ethm
\bpf Let $f_t$ be a deformation of the $n$-ary multiplicative Hom-Nambu-Lie superalgebra $(\g, [\cdot,\cdots,\cdot],\alpha)$ with $f_t=f_{0}+f_rt^r+f_{r+1}t^{r+1}+\cdots,$ i.e., $f_1=f_2=\cdots=f_{r-1}=0.$ Then set $l=r$ in (c$^{'}$), we have
\begin{align*}
&f_{r}(\alpha(x_{1}),\cdots,\alpha(x_{n-1}),[y_{1},\cdots,y_{n}])+[\alpha(x_{1}),\cdots,\alpha(x_{n-1}),f_{r}(y_{1},\cdots,y_{n})]\\
-&\sum_{i=1}^n(-1)^{(|x_{1}|+\cdots+|x_{n-1}|)(|y_{1}|+\cdots+|y_{i-1}|)}\\
\cdot&\Big([\alpha(y_{1}),\cdots,\alpha(y_{i-1}),f_{r}(x_{1},\cdots,x_{n-1},y_{i}),\alpha(y_{i+1}),\cdots,\alpha(y_{n})]\\
+&f_{r}(\alpha(y_{1}),\cdots,\alpha(y_{i-1}),[x_{1},\cdots,x_{n-1},y_{i}],\alpha(y_{i+1}),\cdots,\alpha(y_{n}))\Big)\\
=&0,
\end{align*}
i.e, $\delta^{2}f_{r}(x_{1},\cdots,x_{n-1},[y_{1},\cdots,y_{n}])=0,$
$\delta^{2}(f_r)=0,$ that is, $f_r\in Z^2(\g, \g)_{\bar{0}}.$ By our assumption $H^2(\g, \g)=0$, one gets $f_r\in B^2(\g, \g)_{\bar{0}},$ thus we can find $h_r\in C^1(\g, \g)_{\bar{0}}$ such that $f_r=\delta^{1}h_r.$ Put $\Phi_{t}=\mathrm{id}_{\g}-h_rt^r,$ then $\Phi_{t}\circ(\mathrm{id}_{\g}+h_rt^r+{h_{r}}^{2}t^{2r}+{h_{r}}^{3}t^{3r}+\cdots)
=(\mathrm{id}_{\g}-h_rt^r)\circ(\mathrm{id}_{\g}+h_rt^r+{h_{r}}^{2}t^{2r}+{h_{r}}^{3}t^{3r}+\cdots)=
(\mathrm{id}_{\g}+h_rt^r+{h_{r}}^{2}t^{2r}+{h_{r}}^{3}t^{3r}+\cdots)-(h_rt^r+{h_{r}}^{2}t^{2r}+{h_{r}}^{3}t^{3r}+\cdots)=\mathrm{id}_{\g},$
moreover, $(\mathrm{id}_{\g}+h_rt^r+{h_{r}}^{2}t^{2r}+{h_{r}}^{3}t^{3r}+\cdots)\circ\Phi_{t}=\mathrm{id}_{\g}.$  Hence $\Phi_{t}: \g\rightarrow \g$ is a linear isomorphism and $\Phi_{t}\circ\alpha=\alpha\circ\Phi_{t}.$ Set $f^{'}_{t}(x_{1},\cdots,x_{n})=\Phi_{t}^{-1}f_{t}(\Phi_{t}(x_{1}),\cdots,\Phi_{t}(x_{n})),$ then $f^{'}_{t}$ is also a  deformation of $(\g, [\cdot,\cdots,\cdot],\alpha)$ and $f_{t}\sim f^{'}_{t}.$
Note that $\Phi_{t}f^{'}_{t}(x_{1},\cdots,x_{n})=f_{t}(\Phi_{t}(x_{1}),\cdots,\Phi_{t}(x_{n})).$ Let $f^{'}_{t}=\sum_{i\geq 0}f^{'}_{i}t^{i}.$ Then
$$(\mathrm{id}_{\g}-h_rt^r)\sum_{i\geq 0}f^{'}_{i}(x_{1},\cdots,x_{n})t^{i}=(f_{0}+\sum_{i\geq r}f_{i}t^{i})(x_{1}-h_{r}(x_{1})t^{r},\cdots,x_{n}-h_{r}(x_{n})t^{r}).$$
So
\begin{align*}
&\sum_{i\geq 0}f^{'}_{i}(x_{1},\cdots,x_{n})t^{i}-\sum_{i\geq 0}h_{r}\circ f^{'}_{i}(x_{1},\cdots,x_{n})t^{i+r}\\
=&f_{0}(x_{1},\cdots,x_{n})-\sum_{i=1}^{n}f_{0}(x_{1},\cdots,h_{r}(x_{i}),\cdots,x_{n})t^{r}\\
+&\sum_{1\leq i<j\leq n}f_{0}(x_{1},\cdots,h_{r}(x_{i}),\cdots,h_{r}(x_{j}),\cdots,x_{n})t^{2r}\\
-&\sum_{1\leq i<j<k\leq n}f_{0}(x_{1},\cdots,h_{r}(x_{i}),\cdots,h_{r}(x_{j}),\cdots,h_{r}(x_{k}),\cdots,x_{n})t^{3r}+\cdots\\
+&(-1)^{n}f_{0}(h_{r}(x_{1}),h_{r}(x_{2}),\cdots,h_{r}(x_{n}))t^{nr}+\sum_{i\geq r}f_{i}(x_{1},\cdots,x_{n})t^{i}\\
-&\sum_{i\geq r}\sum_{j=1}^{n} f_{i}(x_{1},\cdots,h_{r}(x_{j}),\cdots,x_{n})t^{i+r}\\
+&\sum_{i\geq r}\sum_{1\leq j\leq k\leq n} f_{i}(x_{1},\cdots,h_{r}(x_{j}),\cdots,h_{r}(x_{k}),\cdots,x_{n})t^{i+2r}+\cdots.
\end{align*}
By the above equation, one gets
$$f^{'}_{0}(x_{1},\cdots,x_{n})=f_{0}(x_{1},\cdots,x_{n})=[x_{1},\cdots,x_{n}];$$
$$f^{'}_{1}(x_{1},\cdots,x_{n})=\cdots=f^{'}_{r-1}(x_{1},\cdots,x_{n})=0;$$
$$f^{'}_{r}(x_{1},\cdots,x_{n})-h_{r}[x_{1},\cdots,x_{n}]=-\sum_{i=1}^{n}[x_{1},\cdots,h_{r}(x_{i}),\cdots,x_{n}]+f_{r}(x_{1},\cdots,x_{n}).$$
Furthermore, we have $$f^{'}_{r}(x_{1},\cdots,x_{n})=-\delta^{1}h_{r}(x_{1},\cdots,x_{n})
+f_{r}(x_{1},\cdots,x_{n})=0,$$
hence, $f^{'}_{t}=f_{0}+\sum_{i\geq r+1}f^{'}_{i}t^{i}.$ By induction, one can prove $f_t\sim f_0,$ that is, $(\g, [\cdot,\cdots,\cdot],\alpha)$ is analytically rigid.
\epf

\end{document}